\documentclass[12pt]{article}

\usepackage{latexsym,amsfonts,amsmath,amssymb,amsthm,url}
\usepackage[english]{babel}
\usepackage[latin1]{inputenc}
\usepackage[T1]{fontenc}
\usepackage{graphics}
\usepackage{color}

\begin{document}

\newtheorem{duge}{Lemma}[section]
\newtheorem{rem}[duge]{Remark}
\newtheorem{prop}[duge]{Proposition}
\newtheorem{defi}[duge]{Definition}
\newtheorem{theo}[duge]{Theorem}
\newtheorem{nota}[duge]{Notation}
\newtheorem{cor}[duge]{Corollary}
\newtheorem{rappel}[duge]{Rappel}
\newtheorem{hypo}[duge]{Hypothesis}
\newtheorem*{ack}{Acknowledgments}

\newcommand{\p}{\bar{p}}
\newcommand{\q}{\bar{q}}
\renewcommand{\P}{\mathbf{P}}
\newcommand{\E}{\mathbf{E}}
\newcommand{\N}{\mathbb{N}}
\newcommand{\Z}{\mathbb{Z}}
\newcommand{\R}{\mathbb{R}}
\newcommand{\B}{B}
\newcommand{\cste}[1]{\ensuremath{c_{#1}}}
\newcommand{\un}{{\mathchoice {\rm 1\mskip-4mu l} {\rm 1\mskip-4mu l} {\rm 1\mskip-4.5mu l} {\rm 1\mskip-5mu l}}}
\newcommand{\defeq}{\ensuremath{\overset{\hbox{\tiny{def}}}{=}}}

\title{\textbf{On the speed of a cookie random walk}}

\author{
\normalsize{\textsc{Anne-Laure Basdevant} and \textsc{Arvind Singh}
\footnote{Address for both authors: \it{Laboratoire de Probabilités
et Modèles Aléatoires, Université Pierre et Marie Curie, 175 rue du
Chevaleret, 75013 Paris, France.} }}}
\date{}
\maketitle \vspace{-1cm}
\begin{center}
University Paris VI
\end{center}

\vspace*{0.2cm}

\begin{abstract}
We consider the model of the one-dimensional cookie random walk when
the initial cookie distribution is spatially uniform and the number
of cookies per site is finite. We give a criterion to decide whether
the limiting speed of the walk is non-zero. In particular, we show that
a positive speed may be obtained for just $3$ cookies per site.
We also prove a result on the continuity of the speed with respect
to the initial cookie distribution.
\end{abstract}

\bigskip
{\small{

\bigskip
\noindent{\bf A.M.S. Classification. }60K35, 60J80, 60F15

\bigskip
\noindent{\bf e-mail. } anne-laure.basdevant@ens.fr,
arvind.singh@ens.fr }}
\section{Introduction}

We consider the model of the multi-excited random walk, also called
cookie random walk, introduced by Zerner in \cite{Zerner05} as a
generalization of the model of the excited random walk described by
Benjamini and Wilson in \cite{Benjamini03} (see also Davis
\cite{Davis99} for a continuous time analogue). The aim of this
paper is to study under which conditions the speed of a cookie
random walk is strictly positive. In dimension $d\ge 2$, this
problem was solved by Kozma \cite{Kozma03,Kozma05} who proved that
the speed is always non-zero. In the one-dimensional case, the speed
can either be zero or strictly positive. We give here a necessary
and sufficient condition to determine if the walk's speed is
strictly positive when the initial cookie environment is
deterministic, spatially uniform and with a finite number of cookies
per site. Let us start with an informal definition of such a
process:

 Let us put $M\geq 1$ "cookies" at each site of $\Z$
and let us pick $p_1,p_2,\dots,p_M \in [\frac{1}{2},1)$. We say that
$p_i$ represents the "strength" of the $i^{\hbox{\tiny{th}}}$ cookie
at any given site. Then, a cookie random walk $X = (X_n)_{n\geq 0}$
is simply a nearest neighbour random walk, eating the cookies it
finds along its path by behaving in the following way:
\begin{itemize}
\item If $X_n = x$ and there is no remaining cookie at site $x$, then $X$
jumps at time $n+1$ to $x +1$ or $x-1$ with equal probability
$\frac{1}{2}$.
\item If $X_n = x$ and there remain the cookies with strengths
$p_j,p_{j+1},\ldots,p_M$ at this site, then $X$ eats the cookie with
attached strength $p_j$ (which therefore disappears from this site)
and then jumps at time $n+1$ to $x+1$ with probability $p_j$ and to
$x-1$ with probability $1-p_j$.
\end{itemize}
This model is a particular case of self-interacting random walk: the
position of $X$ at time $n+1$ depends not only of its position at
time $n$ but also on the number of previous visits to its present
site. Therefore, $X$ is not a Markov process.

Let us now give a formal description of the general model. We define
the set of cookie environments by $\Omega =
[\frac{1}{2},1]^{\N^*\times\Z}$. Thus, a cookie environment is of
the form $\omega = (\omega(i,x))_{i\geq 1,x\in\Z}$ where
$\omega(i,x)$ represents the strength of the $i^{\hbox{\tiny{th}}}$
cookie at site $x$. Given $x\in \Z$ and $\omega\in\Omega$, a cookie
random walk starting from $x$ in the cookie environment $\omega$ is
a process $(X_n)_{n\geq 0}$ on some probability space
$(\mathbf{\Omega},\mathcal{F},\P_{\omega,x})$ such that:
\begin{equation*}
\left\{
\begin{array}{l}
\P_{\omega,x}\{X_0=z\} =1,\\
\P_{\omega,x}\{|X_{n+1}- X_n|= 1\}=1,\\
\P_{\omega,x}\{X_{n+1} = X_n + 1\hbox{ | }X_1,\ldots,X_n\} =
\omega(j,X_n) \hbox{ where } j = \sharp\{0\leq i\leq n\hbox{ ,
}X_i=X_n\}.
\end{array}
\right.
\end{equation*}
In this paper, we restrict our attention to the set of environments
$\Omega^u_{M}\subset\Omega$ which are spatially uniform with at most
$M\geq 1$ cookies per site:
\begin{equation*}
\omega\in\Omega^u_M\quad\Longleftrightarrow\quad\left\{
\begin{array}{l}
\hbox{for all $x\in\Z$ and all $i\geq 1$ } \omega(i,x)=\omega(i,0),\\
\hbox{for all $i>M$ } \omega(i,0)=\frac{1}{2},\\
\hbox{for all $i\geq 1$ } \omega(i,0)<1.
\end{array}\right.
\end{equation*}
The last condition $\omega(i,0)<1$ is introduced only to exclude
some possible degenerated cases but can be relaxed (see Remark
\ref{RemDegenerate}). A cookie environment $\omega\in\Omega^u_M$ may
be represented by $(M,\p)$ where
\begin{equation*}
\begin{array}{l}
\p = (p_1,\ldots,p_M) = (\omega(1,0),\ldots,\omega(M,0)).
\end{array}
\end{equation*}
In this case, we shall say that the associated cookie random walk is
an $(M,\p)$-cookie random walk and we will use the notation
$\P_{(M,\p)}$ instead of $\P_{\omega}$.

The question of the recurrence or transience of a cookie random walk
was solved by Zerner in \cite{Zerner05} for general cookie
environments (even in the case where the initial cookie environment may itself be random).
In particular, he proved that, if $X$ is an $(M,\p)$
cookie random walk, there is a phase transition according to the
value of
\begin{equation}\label{DefAlpha}
\alpha = \alpha(M,\p) \defeq \sum_{i=1}^M (2p_i-1)-1.
\end{equation}
\begin{itemize}
\item If $\alpha \leq 0$ then the walk is recurrent \emph{i.e.} $\limsup
X_n = -\liminf X_n = +\infty$ a.s.
\item If $\alpha > 0$ then $X$ is transient toward $+\infty$ \emph{i.e}
$\lim X_n = +\infty$ a.s.
\end{itemize}
In particular, for $M=1$, the cookie random walk is always recurrent
for any choice of $\p$. However, as soon as $M\geq 2$, the cookie
random walk can either be transient or recurrent depending on $\p$.
Zerner \cite{Zerner05} also proved that the speed of a
$(M,\p)$-cookie random walk $X$ is always well defined (but may or
not be zero). Precisely,
\begin{itemize}
\item there exists a constant $v(M,\p) \geq 0$ such that
\begin{equation*}
\frac{X_n}{n}\underset{n\to\infty}{\longrightarrow}v(M,\p)\quad\hbox{$\P_{(M,\p)}$-almost
surely.}
\end{equation*}
\item The speed is monotonic in $\p$: if $\p = (p_1,\ldots,p_M)$ and $\q =
(q_1,\ldots,q_M)$ are two cookie environments such that $p_i \leq
q_i$ for all $i$, then $v(M,\p) \leq v(M,\q)$.
\item The speed of a $(2,\p)$-cookie random walk is always $0$.
\end{itemize}

The question of whether one can construct a $(M,\p)$-cookie random
walk with strictly positive speed was affirmatively answered by
Mountford, Pimentel and Valle \cite{Mountford06} who considered the
case where all the cookies have the same strength
$p\in[\frac{1}{2},1)$ \emph{i.e.} the cookie vector $\p$ has the
form $[p]_M \defeq (p,\ldots,p)$. They showed that:
\begin{itemize}
\item For any $p\in (\frac{1}{2},1)$, there exists an $M_0$ such that for
all $M>M_0$ the speed of the $(M,[p]_M)$-cookie random walk is
strictly positive.
\item If $M(2p-1)<2$, then the speed of the $(M,[p]_M)$-cookie random walk is zero.
\end{itemize}
They also conjectured that when $M(2p-1)> 2$, the speed should be
non-zero. The aim of this paper is to prove that such is indeed the
case.
\begin{theo}\label{MainTheo1} Let $X$ denote a $(M,\p)$-cookie random walk, then
\begin{equation*}
\lim_{n\to\infty}\frac{X_n}{n} = v(M,\p) >0 \quad
\Longleftrightarrow \quad \alpha(M,\p) > 1
\end{equation*}
where $\alpha(M,\p)$ is given by (\ref{DefAlpha}).
\end{theo}
In particular, we see that a non-zero speed may be achieved for as
few as $3$ cookies per site. Comparing this result with the
transience/recurrence criteria, we have a second order phase
transition at the critical value $\alpha=1$. In fact, it shall be proved
in a forthcoming paper that, for $0<\alpha<1$, the rate of
transience of $X_n$ is of order  $n^{\frac{\alpha+1}{2}}$.

One would certainly like an explicit calculation of the limiting
velocity in term of the cookie environment $(M,\p)$ but this seems a
challenging problem (one can still look at the end of Section $3$
where we give an implicit formula for the speed). However, one can
prove that the speed is continuous in $\p$ and has a positive right
derivative at all its critical points:
\begin{theo}\label{MainTheo2}
\begin{itemize}
\item For each $M$, the speed $v(M,\p)$ is a continuous function of $\p$ in $\Omega^u_M$.
\item For any environment $(M,\p_c)$ with $\alpha(M,\p_c) = 1$,
there exists a constant $C>0$ (depending on $(M,\p_c)$) such that
\begin{equation*}
\lim_{
\begin{subarray}{l}
\p\to\p_c\\
\p\in\Omega^u_M\\
\alpha(\p)>1
\end{subarray}
}\frac{v(M,\p)}{\alpha(M,\p)-1} = C.
\end{equation*}
\end{itemize}
\end{theo}
In particular, for $M\geq 3$, the (unique) critical value for an
$(M,[p]_M)$-cookie random walk is $p_c = \frac{1}{M} + \frac{1}{2}$
and the function $v(p)$ is continuous, non-decreasing, zero for
$p\leq p_c$, and admits a finite strictly positive right derivative
at $p_c$.

The remainder of this paper is organized as follow. In the next
section, we construct a Markov process associated with the hitting
time of the cookie random walk. The method is similar to that used
by Kesten, Kozlov and Spitzer \cite{Kesten75} for the determination
of the rates of transience of a random walk in a one-dimensional
random environment. It turns out that, in our setting, the resulting
process is a branching process with random migration. The study of
this process and of its stationary distribution is done in Section
$3$. This enables us to complete the proof of Theorem
\ref{MainTheo1}. Finally, the last section is dedicated to the proof
of Theorem \ref{MainTheo2}.

\section{An associated branching process with migration}\label{defZn}
In the remainder of this paper, $X=(X_n)_{n\geq 0}$ will denote a
$(M,\p)$-cookie random walk. Since the speed of a recurrent cookie
random walk is zero, we will also assume that we are in the
transient regime \emph{i.e.}
\begin{equation}
\alpha(M,\p) =\sum_{i=1}^{M}(2p_i-1) -1 > 0.
\end{equation}
For the sake of brevity, we simply write $\P_{x}$ for
$\P_{(M,\p),x}$ and $\P$ instead of $\P_{0}$ (the process starting
from $0$). Let $T_n$ stand for the hitting time of level $n\geq 0$
by $X$:
\begin{equation}
T_n= \inf(k\geq 0\hbox{ , }X_k = n).
\end{equation}
For $0\leq k\leq n$, let $U_i^n$ denote the number of jumps of the
cookie random walk from site $i$ to site $i-1$ before reaching level
$n$
\begin{equation*}
U_i^n=\sharp\{0\leq k < T_n , \;X_{k}=i \mbox{ and } X_{k+1}=i-1\}.
\end{equation*}
Let also $K_n$ stand for the total time spent by $X$ in the negative
half-line up to time $T_n$
\begin{equation*}
K_n= \sharp\{0\leq k \leq T_n, \; X_k <0\}.
\end{equation*}
A simple combinatorial argument readily yields
\begin{equation*}
T_n = K_n-U_0^n + n + 2\sum_{k=0}^{n}U^n_k.
\end{equation*}
Notice that, as $n$ tends to infinity, the random variable $K_n$
increases almost surely toward $K_\infty$, the total time spent by
the cookie random walk in the negative half line. Similarly, $U_0^n$
increases toward $U_0^\infty$ the total number of jumps from $0$ to
$-1$. Since $X$ is transient, $K_\infty+U_0^\infty$ is almost-surely
finite and therefore
\begin{equation}\label{equivTU}
T_n \underset{n\to\infty}{\sim} n + 2\sum_{k=0}^{n}U^n_k.
\end{equation}
Let us now prove that for each $n$, the reverse process
$(U^n_n,U^n_{n-1},\ldots, U^n_1,U^n_0)$ has the same law as the $n$
first steps of some branching process $Z$ with random migration. We
first need to introduce some notations. Let $(B_i)_{i\geq 1}$ denote
a sequence of independent Bernoulli random variable under $\P$ with
distribution:
\begin{equation}\label{DefBernoulli}
\P\{B_i = 1\} = 1-\P\{B_i = 0\} = \left\{
\begin{array}{ll}
p_i&\hbox{ if $i\leq M$,}\\
\frac{1}{2}&\hbox{ if $i> M$.}
\end{array}\right.
\end{equation}
For $j\in\N$, define
\begin{equation*}
k_j = \min(k\geq 1 , \sharp\{1\leq i \leq k, B_i =1\}=j+1)
\end{equation*} and
\begin{equation*}
A_j = \sharp\{1\leq i \leq k_j, B_i =0\} = k_j - j - 1.
\end{equation*}
We have the following easy lemma:
\begin{duge}\label{LemmaA}\begin{itemize}
\item For any $i,j\geq 0$, we have $\P\{A_j=i\}>0$.
\item For all $j\geq M$, we have
\begin{equation}
A_j \overset{\hbox{\tiny{law}}}{=} A_{M-1} + \xi_1 + \ldots +
\xi_{j-M+1}
\end{equation}
where $(\xi_i)_{i\geq 0}$ is a sequence of i.i.d. geometrical random
variable with parameter $\frac{1}{2}$ independent of $A_{M-1}$.
\end{itemize}
\end{duge}
\begin{proof}
The first part of the lemma is a direct consequence of the
assumption that $\p$ is such that $p_k<1$ for all $k$. To prove the
second part, we simply notice that $k_{M-1}\geq M$ so that for
$j\geq M$, the random variable $A_{j} - A_{M-1}$ has the same law as
the random variable
\begin{equation}\label{souduge1}
\min(k\geq 1 , \sharp\{1\leq i \leq k, \widetilde{B}_i =1\}=j+1-M) -
j-1 + M
\end{equation}
where $(\widetilde{B}_i)_{i\geq 0}$ is a sequence of i.i.d. random
variables independent of $A_{M-1}$ and with common Bernoulli
distribution
$\P\{\widetilde{B}_i=0\}=\P\{\widetilde{B}_i=1\}=\frac{1}{2}$. It is
clear that (\ref{souduge1}) has the same law as $\xi_1 + \ldots +
\xi_{j-M+1}$.
\end{proof}

By possibly extending the probability space, we now construct a
process $Z=(Z_n,n\ge 0)$ and a family of probability
$(\mathbb{P}_z)_{z\geq 0}$ such that, under $\mathbb{P}_z$, the
process $Z$ is a Markov chain starting from $z$, with transition
probability:
\begin{equation*} \left\{
\begin{array}{l}
\mathbb{P}_z\{Z_0=z\}=1,\\
\mathbb{P}_z\{Z_{n+1}=k \;|\;Z_n=j\}=\P\{A_{j}=k\}.
\end{array}\right.
\end{equation*}
Since the family of probabilities $(\mathbb{P}_z)$ depends on the
law of the cookie environment $(M,\p)$, we should rigourously write
$\mathbb{P}_{(M,\p),z}$ instead of $\mathbb{P}_z$. However, when
there is no possible confusion we will keep using the abbreviated
notation. Furthermore,  we will simply write $\mathbb{P}$ for
$\mathbb{P}_0$ and $\mathbb{E}$ will stand for the expectation with
respect to $\mathbb{P}$.

Let us now notice that, in view of the previous lemma, $Z_n$ under
$\mathbb{P}_z$ may be interpreted as the number of particles alive
at time $n$ of a branching process with random migration starting
from $z$, that is a branching process which allows immigration and
emigration (see Vatutin and Zubkov \cite{Vatutin93} for a survey on
these processes). Indeed:
\begin{itemize}
\item If $Z_n = j \geq M-1$, then according to Lemma \ref{LemmaA}, $Z_{n+1}$ has the same law as
$\sum_{k=1}^{j - M +1}\xi_k + A_{M-1}$, \emph{i.e.} $M-1$ particles
emigrate and the remaining particles reproduce according to a
geometrical law with parameter $\frac{1}{2}$ and there is also an
immigration of $A_{M-1}$ new particles.
\item If $Z_n = j \in \{0,\ldots,M-2\}$ then $Z_{n+1}$ has the same law as
$A_{j}$ \emph{i.e.} all the $j$ particles emigrate and $A_j$ new
particles immigrate.
\end{itemize}
We can now state the main result of this section:
\begin{prop}\label{propUegalZ}
For each $n\in \N$, $(U_n^n,U_{n-1}^n,\ldots,U_0^n)$ under $\P$ has
the same law as $(Z_0,Z_1,\ldots,Z_n)$ under $\mathbb{P}$.
\end{prop}
\begin{proof}
The argument is similar to the one given by Kesten \emph{et al.} in
\cite{Kesten75}. Recall that $U_i^n$ represents the numbers of jumps
of the cookie random walk $X$ from $i$ to $i-1$ before reaching $n$.
Then, conditionally on $(U_n^n,U_{n-1}^n,\ldots,U_{i+1}^n)$, the
number of jumps $U^n_i$ from $i$ to $i-1$ depends only on the number
of jumps from $i+1$ to $i$, that is, depends  only of $U^n_{i+1}$.
This shows that $(U_n^n,U_{n-1}^n,\ldots,U_0^n)$ is indeed a Markov
process.

By definition, $Z_0=0$ $\mathbb{P}$-a.s. and $U_n^n =0$ $\P$-a.s. It
remains to compute $\P\{U_i^n=k\;|\;U_{i+1}^n=j\}$. Note that the
number of jumps from $i$ to $i-1$ before reaching level $n$ is equal
to the number of jumps from $i$ to $i-1$ before reaching $i+1$ for
the first time plus the sum of the number of jumps from $i$ to $i-1$
between two consecutive jumps from $i+1$ to $i$ which occur before
reaching level $n$. Thus, conditionally on $\{U_{i+1}^n=j\}$, the
random variable $U^n_{i}$ has the same law as the number of failures
(i.e. $B_k=0$) in the Bernoulli sequence $(B_1,B_2,B_3,\ldots)$
defined by (\ref{DefBernoulli}) before having exactly $j+1$
successes. This is precisely the definition of $A_j$ and therefore
$\P\{U_i^n=k\;|\;U_{i+1}^n=j\} = \mathbb{P}_j\{Z_1=k\}$.
\end{proof}

Since $U^n_0$ is the number of jumps from $0$ to $-1$ of the cookie
random walk $X$ before reaching level $n$ and since we assumed that
the cookie random walk $X$ is transient, $U_0^n$ increases almost
surely toward the total number $U^\infty_0$ of jumps of $X$ from $0$
to $-1$. In view of the previous proposition, this implies that
under $\mathbb{P}$, $Z_n$ converges in law toward a random variable
which we denote by $Z_\infty$.

Let us also note that $Z$ is a irreducible Markov chain (this is a
consequence of part $1$ of Lemma \ref{LemmaA}). Since $Z$ converges
in law toward a limiting distribution, this shows that $Z$ is in
fact a positive recurrent Markov chain. In particular, $Z_n$
converges in law toward $Z_\infty$ independently of its starting
point (\emph{i.e.} the law of $Z_\infty$ is the same under any
$\mathbb{P}_x$) and the law of $Z_\infty$ is also the unique
invariant probability for $Z$.

\begin{cor}\label{CorVitesse}
Recall that $v(M,\p)$ denotes the limiting speed of the cookie
random walk $X$. We have
\begin{equation*}
v(M,\p) = \frac{1}{1+2\mathbb{E}[Z_\infty]}\quad\hbox{(with the
convention $0=\frac{1}{+\infty}$).}
\end{equation*}
In particular, the speed of an $(M,\p)$-cookie random walk is non
zero i.i.f. the limiting random variable $Z_\infty$ of its
associated process $Z$ has a finite expectation.
\end{cor}

\begin{proof}
Since $X$ is transient, we have the well known equivalence valid for
$v\in [0,\infty]$ :
\begin{equation}
\frac{X_n}{n}\underset{n\to\infty}{\longrightarrow}v\quad\hbox{$\P$-a.s.}
\qquad \Longleftrightarrow \qquad
\frac{T_n}{n}\underset{n\to\infty}{\longrightarrow}\frac{1}{v}\quad\hbox{$\P$-a.s.}
\end{equation}
On the one hand, this equivalence and (\ref{equivTU}) yield
\begin{equation}\label{this1}
\frac{1}{n}\sum_{k=0}^n
U^n_k\underset{n\to\infty}{\longrightarrow}\frac{1}{2v(M,\p)} -
\frac{1}{2}\quad\hbox{$\P$-a.s.}
\end{equation}
On the other hand, making use of an ergodic theorem for the positive
recurrent Markov chains $Z$ with stationary limiting distribution
$Z_\infty$, we find that
\begin{equation}\label{this2}
\frac{1}{n}\sum_{i=1}^{n}Z_k\underset{n\to\infty}{\rightarrow}\mathbb{E}[Z_\infty]\quad\hbox{$\mathbb{P}$-a.s.}
\end{equation}
(this result is valid even if $\mathbb{E}[Z_\infty]=\infty$).
Proposition \ref{propUegalZ} implies that the limits in
(\ref{this1}) and (\ref{this2}) are the same. This completes the
proof of the corollary.
\end{proof}

\begin{rem}\label{RemDegenerate} We assumed in the definition of an
$(M,\p)$ cookie environment that
\begin{equation*}
p_i \neq 1 \quad\hbox{ for all $1\leq i\leq M$.}
\end{equation*}
This hypothesis is intended only to ensure that $Z$ starting from
$0$ is not almost surely bounded (for instance, if $p_1=1$ then $0$
is a absorbing state for $Z$). More generally, one may check from
the definition of the random variables $A_j$ that $Z$ starting from
$0$ is almost surely unbounded i.i.f.
\begin{equation}\label{condiIrreductible}
\sharp\{1\leq j \leq i \hbox{ , } p_j = 1 \} \leq \frac{i}{2} \quad \hbox{for all $1\leq i\leq M$.}
\end{equation}
When this condition fails, $Z$ starting from $0$ is almost surely bounded by $M-1$, thus $\mathbb{E}[Z_\infty]<\infty$ and the speed of the associated cookie random walk is strictly positive. Otherwise, when (\ref{condiIrreductible}) is fulfilled, $Z$ ultimately hits any level $x\in\N$ with probability $1$ and the proof of Theorem \ref{MainTheo1} remains valid.
\end{rem}

\section{Study of $Z_\infty$.}
We proved in the previous section that the strict positivity of the
speed of the cookie random walk $X$ is equivalent to the existence
of a finite first moment for the limiting distribution of its
associated Markov chain $Z$. We shall now show that, for any cookie
environment $(M,\p)$ (with $\alpha(M,\p)>0$), we have
\begin{equation*}
\mathbb{E}[Z_\infty] \defeq \mathbb{E}_{(M,\p)}[Z_\infty]<\infty
\qquad \Longleftrightarrow \qquad \alpha(M,\p)>1.
\end{equation*}
This will complete the proof of Theorem \ref{MainTheo1}. We start by
proving that $Z_\infty$ cannot have moments of any order.
\begin{prop}\label{espZM} We have
\begin{equation*}
\mathbb{E}\left[Z_\infty^{M-1}\right] = +\infty.
\end{equation*}
\end{prop}

\begin{proof}  Let us introduce the first return time to $0$ for
$Z$:
\begin{equation*}
\sigma = \inf(n\geq 1\; ,\;Z_n=0).
\end{equation*}
Since $Z$ is a positive recurrent Markov chain, we have
$1\leq\mathbb{E}_0[\sigma] < \infty$ and the invariant probability
measure is given for any $y\in \N$ by
$$\mathbb{P}\{Z_\infty=y\}=\frac{\mathbb{E}_0\left[\sum_{k=0}^{\sigma-1}\un_{Z_k=y}\right]}{\mathbb{E}_0[\sigma]}.$$
 A monotone convergence argument yields
\begin{equation}\label{PropCou1}
\mathbb{E}_{0}\left[\sum_{k=0}^{\sigma-1}Z_k^{M-1}\right] =
\mathbb{E}_{0}[\sigma]\mathbb{E}[Z_\infty^{M-1}]
\end{equation}
(where both side of this equality may be infinite). We can find
$n_0\in\N^*$ such that $\mathbb{P}_{0}\{Z_{n_0} = M,\; n_0 <
\sigma\}> 0$ (in fact, since we assume that $p_i<1$ for all $i$, we
can choose $n_0=1$). Therefore, making use of the Markov property of
$Z$, we find that
\begin{eqnarray}
\nonumber\mathbb{E}_{0}\left[\sum_{k=0}^{\sigma-1}Z_k^{M-1}\right]
&\geq& \mathbb{P}_{0}\{Z_{n_0} = M,\; n_0<\sigma\}\mathbb{E}_{M}\left[\sum_{k=0}^{\sigma-1}Z_k^{M-1}\right]\\
\label{PropCou2} &=& \mathbb{P}_{0}\{Z_{n_0} = M,\; n_0 < \sigma
\}\sum_{k=0}^{\infty}\mathbb{E}_{M}\left[Z_{k\wedge\sigma}^{M-1}\right].
\end{eqnarray}
In view of (\ref{PropCou1}) and (\ref{PropCou2}), we just need to
prove that
\begin{equation}\label{PropCou3}
\sum_{k=0}^{\infty}\mathbb{E}_{M}\left[Z_{k\wedge\sigma}^{M-1}\right]
= \infty.
\end{equation}
We now use a coupling argument. Let us define a new Markov chain
$\widetilde{Z}$ such that, under $\mathbb{P}_z$, the process evolves
in the following way
\begin{itemize}
\item $\widetilde{Z}_0= z$,
\item if $\widetilde{Z}_n = k \in \{0,1,\ldots,M-1\}$ then
$\widetilde{Z}_{n+1}=0$,
\item if $\widetilde{Z}_n = k > M-1$ then $\widetilde{Z}_{n+1}$ has
the same law as $\sum_{i=1}^{k-(M-1)}\xi_i$ where $(\xi_i)_{i\ge 1}$
is a sequence of i.i.d. geometrical random variables with parameter
$\frac{1}{2}$.
\end{itemize}
Thus, $\widetilde{Z}$ is a branching process with emigration: at
each time $n$, there are $\min(\widetilde{Z}_n,M-1)$ particles which
emigrate the system and the remaining particles reproduce according
to a geometrical law of parameter $\frac{1}{2}$.

Recall that $Z$ is a branching process with migration, where at most
$M-1$ particles emigrate at each unit of time, and has the same
offspring reproduction law as $\widetilde{Z}$. Therefore, for any
$z\geq 0$, under $\mathbb{P}_z$, the process $\widetilde{Z}$ is
stochastically dominated by $Z$. Since $0$ is an absorbing state for
$\widetilde{Z}$, this implies that, for all $n\geq 0$ and all $z\geq
0$,
\begin{equation}\label{PropCou4}
\mathbb{E}_{z}[\widetilde{Z}^{M-1}_n] \leq
\mathbb{E}_{z}[Z^{M-1}_{n\wedge\sigma}].
\end{equation}
Our process $\widetilde{Z}$ belongs to the class of processes
studied by Kaverin \cite{Kaverin88}. Moreover,  all the hypotheses
of Theorem $1$ of \cite{Kaverin88} are clearly fulfilled (in the
notation of \cite{Kaverin88}, we have here $\lambda=\theta=M-1$ and
$B=1$). Therefore, for any $z\geq M$, there exists a constant $c >0$
(depending on $z$) such that
\begin{equation}\label{PropCou5}
\mathbb{E}_{z}[\widetilde{Z}^{M-1}_n]\underset{n\to\infty}{\sim}\frac{c}{n}.
\end{equation}
The combination of (\ref{PropCou4}) and (\ref{PropCou5}) yield
(\ref{PropCou3}).
\end{proof}

\begin{rem} In view of the last proposition and Corollary \ref{CorVitesse}, we recover the fact
that for $M=2$, the speed of the cookie random walk is always zero.
\end{rem}

In order to study more precisely the distribution of $Z_\infty$, we
will need the following lemma
\begin{duge}\label{LemmeAM} We have
\begin{equation*}
\E\left[A_{M-1}\right] = 2\sum_{i=1}^{M}(1-p_i).
\end{equation*}
\end{duge}

\begin{proof} Recall that $(B_i)_{i\geq 1}$ denotes a sequence of independent Bernoulli
random variables with distribution given by (\ref{DefBernoulli}).
Recall also that
\begin{eqnarray*}
k_{M-1} &=& \min(k\geq 1 , \sharp\{1\leq i \leq k, B_i =1\}=M),\\
A_{M-1} &=& k_{M-1} - M.
\end{eqnarray*}
In particular, for any $j\in \N^*$
\begin{eqnarray*}
\P\{A_{M-1}=j\}&=&\P\{k_{M-1}=M+j\}\\
&=&\P\Big\{\sharp\{1\leq i \leq M+j-1, B_i
=1\}=M-1\Big\}\P\{B_{M+j}=1\}.
\end{eqnarray*}
 Hence,
\begin{multline*}
\P\{A_{M-1}=j\}\\
\begin{aligned} &=\!
\frac{1}{2}\!\sum_{l=1}^{M}\P\Big\{\sharp\{1\!\leq\!i\!\leq\!M, B_i
=1\}=M\!-\!l\Big\} \P\!\Big\{\sharp\{M\!+\!1\leq i \leq M\!+\!j\!-\!1, B_i =1\}=l\!-\!1\Big\}\\
&=\sum_{l=1}^{M\wedge j}\P\Big\{\sharp\{1\leq i \leq M, B_i
=1\}=M-l\Big\} C_{j-1}^{l-1}\left(\frac{1}{2}\right)^j.
\end{aligned}
\end{multline*}
Let $L=\sharp\{1\leq i \leq M, B_i =0\}$, therefore
\begin{equation*}
\P\{A_{M-1}=j\}=\sum_{l=1}^{M\wedge j}\P\{L=l\}
C_{j-1}^{l-1}\left(\frac{1}{2}\right)^j \quad \mbox{ for } j\in
\N^*.
\end{equation*}
Making use of the relation $j C_{j-1}^{l-1}=lC_j^l$, we get
\begin{eqnarray*}
\E[A_{M-1}]&=&\sum_{j=1}^{\infty}\sum_{l=1}^{M\wedge j} \P\{L=l\}\left(\frac{1}{2}\right)^j l C_{j}^{l}\\
&=&\sum_{l=0}^{M}l \P\{L=l\}\sum_{j=l}^{\infty}\left(\frac{1}{2}\right)^j  C_{j}^{l}\\
&=& 2 \sum_{l=0}^{M}l \P\{L=l\}\\
&=&2\E[L].
\end{eqnarray*}
We now compute $\E[L]$ by induction on the number of cookies.
$$\P\{L=l\}=\sum_{1\le i_1<i_2<\ldots<i_l\le M} \prod_{j=1}^l(1-p_{i_j})
\prod_{j\notin\{ i_1,\ldots,i_l\}}p_j \qquad \mbox{ for } 0\le l\le
M.$$ Decomposing the last sum according to whether $i_1=1$ or
$i_1\neq1$, we obtain
\begin{multline*}
\E[L]\\
\begin{aligned}
&= \!\sum_{l=1}^M l\hspace*{-0.2cm} \sum_{ 2\le
i_2<\ldots<i_l\le M}\hspace*{-0.5cm} (1-p_1)\prod_{j=2}^l(1-p_{i_j})
\hspace*{-0.5cm}\prod_{j\notin\{
1,i_2,\ldots,i_l\}}\hspace*{-0.5cm}p_j +\sum_{l=0}^{M-1}
l\hspace*{-0.2cm}\sum_{2\le i_1<\ldots<i_l\le M}
\prod_{j=1}^l(1-p_{i_j}) p_1 \hspace*{-0.3cm}\prod_{j\notin\{
1,i_1,\ldots,i_l\}}\hspace*{-0.5cm}p_j\\
&= \sum_{l=0}^{M-1} (l-(1-p_1))\hspace*{-0.2cm} \sum_{ 2\le
i_1<\ldots<i_l\le M} \prod_{j=1}^l(1-p_{i_j})
\hspace*{-0.5cm}\prod_{j\notin\{
1,i_1,\ldots,i_l\}}\hspace*{-0.5cm}p_j\\
&= \sum_{l=0}^{M-1} (l-(1-p_1))\P\{\tilde{L}=l\}\\
&= \E[\tilde{L}]+1-p_1,
\end{aligned}
\end{multline*}
where $\tilde{L}=\sharp\{2\leq i \leq M, B_i =0\}$. Finally, we
conclude by induction that
$$\E[L]= \sum_{i=1}^{M}(1-p_i).$$
\end{proof}

We now study the law of the limiting distribution $Z_\infty$ of the
Markov chain $Z$. This is done via the study of its probability
generating function (p.g.f.)
\begin{equation*}
G(s) = \mathbb{E}\left[s^{Z_\infty}\right]\quad\hbox{ for
$s\in[0,1]$.}
\end{equation*}

\begin{duge}\label{relationG}
The p.g.f. $G$ of $Z_\infty$ is the unique p.g.f. solution of the
following equation
\begin{equation}\label{equafonc}1-G\left(\frac{1}{2-s}\right)=a(s)(1-G(s))+b(s) \qquad
\mbox{ for all $s\in [0,1]$,} \end{equation} with
$$a(s)=\frac{1}{(2-s)^{M-1}\E\left[s^{A_{M-1}}\right]},$$
and
$$b(s)=1-\frac{1}{(2-s)^{M-1}\E\left[s^{A_{M-1}}\right]}+
\sum_{k=0}^{M-2}G^{(k)}(0)\left(\frac{\E\left[s^{A_{k}}\right]}{(2-s)^{M-1}\E\left[s^{A_{M-1}}\right]}-\frac{1}{(2-s)^k}
 \right).$$
\end{duge}
\begin{proof} The law of $Z_\infty$ is a stationary distribution for
the Markov chain $Z$, therefore
\begin{eqnarray*}
G(s)\;=\;\mathbb{E}\left[\mathbb{E}_{Z_\infty}\left[s^{Z_1}\right]\right] &=&
\sum_{k=0}^{\infty} \mathbb{P}\{Z_\infty=k\}\mathbb{E}_k\left[s^{Z_1}\right]\\
&=& \sum_{k=0}^{M-2}
\mathbb{P}\{Z_\infty=k\}\mathbb{E}_k\left[s^{Z_1}\right] +
\sum_{k=M-1}^{\infty}
\mathbb{P}\{Z_\infty=k\}\mathbb{E}_k\left[s^{Z_1}\right].
\end{eqnarray*}
By definition of $Z$, for $0\leq k\leq M-2$,  $Z_1$ under
$\mathbb{P}_{k}$ has the same law as $A_k$ under $\P$. For $k\geq
M-1$, $Z_1$ under $\mathbb{P}_{k}$ has the same law as
$A_{M-1}+\xi_{1}+\ldots+\xi_{k-M+1}$ where $(\xi_i)_{i\ge 1}$ is a
sequence of i.i.d. random variables independent of $A_{M-1}$ and
with geometric distribution with parameter $\frac{1}{2}$. Thus,
\begin{eqnarray*}
G(s) &=&
\sum_{k=0}^{M-2}\mathbb{P}\{Z_\infty=k\}\E\left[s^{A_{k}}\right]+
\sum_{k=M-1}^{\infty}\mathbb{P}\{Z_\infty=k\}\E\left[s^{A_{M-1}+\xi_1+\ldots+\xi_{k+1-M}}\right]\\
&=&\sum_{k=0}^{M-2}\mathbb{P}\{Z_\infty=k\}\E\left[s^{A_{k}}\right]+
\frac{\E\left[s^{A_{M-1}}\right]}{\E\left[s^{\xi}\right]^{M-1}}
\sum_{k=M-1}^{\infty}\mathbb{P}\{Z_\infty=k\}\E\left[s^{\xi}\right]^{k}\\
&=&\sum_{k=0}^{M-2}\mathbb{P}\{Z_\infty=k\}\left(\E\left[s^{A_{k}}\right]-
\E\left[s^{A_{M-1}}\right]\E\left[s^{\xi}\right]^{k+1-M}
\right)+\frac{\E\left[s^{A_{M-1}}\right]}{\E\left[s^{\xi}\right]^{M-1}}G\left(\E\left[s^{\xi}\right]\right).
\end{eqnarray*}
Since $\E\left[s^{\xi}\right]=\frac{1}{2-s},$ and
$\mathbb{P}\{Z_\infty=k\} = G^{k}(0)$, we get
$$G(s)=\sum_{k=0}^{M-2}G^{k}(0)\left(\E\left[s^{A_{k}}\right]-\E\left[s^{A_{M-1}}\right](2-s)^{M-1-k}
 \right)+\E\left[s^{A_{M-1}}\right](2-s)^{M-1}G\left(\frac{1}{2-s}\right),$$
from which we deduce that $G$ solves (\ref{equafonc}). Furthermore,
the uniqueness of the solution of this equation amongst the class of
probability generating function is a direct consequence of the
uniqueness of the stationary law for the irreducible Markov chain
$Z$.
\end{proof}

Given two functions $f$ and $g$, we use the classical notation $f(x)
= \mathcal{O}(g(x))$ in the neighbourhood of zero if $|f(x)| \leq
C|g(x)|$ for some constant $C$ and all $|x|$ small enough.
\begin{duge}\label{LemmaFormeAB} The functions $a$ and $b$ of Lemma \ref{relationG} are
analytic on $(0,2)$. In particular, they admit a Taylor expansion of
any order near point $1$ and, as $x$ goes to $0$:
\begin{eqnarray*}
a(1-x)&=&1-\alpha x + \mathcal{O}(x^2),\\
b(1-x)&=& \mathcal{O}(x).
\end{eqnarray*}
\end{duge}
\begin{proof} Recall the definitions of the random variables
$A_k$ given in Section $2$. Since a geometric random variable with
parameter $\frac{1}{2}$ admits exponential moments of order strictly
smaller than $2$, it follows that the p.g.f. $s\mapsto\E[s^{A_{k}}]$
are strictly positive and analytic on $(0,2)$. From the explicit
form of the functions $a$ and $b$ given in the previous lemma, we
conclude that these two functions are indeed analytic on $(0,2)$. A
Taylor expansion of $a$ near $1$ gives
\begin{equation}\label{deva}
a(1-x)=1-\left(M-1-\E[A_{M-1}]\right)x+\mathcal{O}(x^2)=1-\alpha x
+\mathcal{O}(x^2),
\end{equation}
where we used Lemma \ref{LemmeAM} for the last equality. Since $G$
is a p.g.f. we have $G(1)=1$ which, in view of (\ref{equafonc}),
yields $b(1)=0$ and therefore $b(1-x)= \mathcal{O}(x)$.
\end{proof}

The following proposition relies on a careful study of equation (\ref{equafonc}) and is the
key to the proof of Theorem \ref{MainTheo1}.
\begin{prop}\label{propG}
Recall that
\begin{equation*}
\alpha = \sum_{i=1}^M (2p_i-1)-1 > 0.
\end{equation*}
The p.g.f. $G$ of $Z_\infty$ is such that, as $x$ goes to $0$:
\begin{itemize}
\item if $0<\alpha<1$,  then $1-G(1-x)\sim \cste{1}x^\alpha$, for some constant $\cste{1}>0$.

In particular $\E[Z_\infty]=+\infty$.
\item if $\alpha=1$,  then $1-G(1-x)\sim \cste{2}x|\ln x|$, for some constant $\cste{2}>0$.

In particular $\E[Z_\infty]=+\infty$.
\item if $\alpha>1$,  then $1-G(1-x)= \cste{3}x+\mathcal{O}(x^{2\wedge
\alpha})$ for some constant $\cste{3}>0$.

In particular $\E[Z_\infty]<+\infty$.
\end{itemize}

\end{prop}

\begin{proof} Since $G$ is a p.g.f, it is completely monotonic and we just need to prove the
proposition along the sequence $x=\frac{1}{n}$ with $n\in\N^*$.
Making use of Lemma \ref{relationG} with $s=1-\frac{1}{n}$, we get,
for all $n\geq 1$
\begin{equation*}
1-G\left(1-\frac{1}{n+1}\right)=a\left(1-\frac{1}{n}\right)
\left(1-G\left(1-\frac{1}{n}\right)\right)+b\left(1-\frac{1}{n}\right).
\end{equation*}
Let us define the sequence $(u_n)_{n\ge 1}$ by
\begin{equation}\label{defun}
\left\{
\begin{array}{l}
u_1=1-G(0)=1-\P(Z_\infty=0)>0,\\
u_n=\frac{1-G(1-1/n)}{\prod_{i=1}^{n-1}a(1-1/i)}\quad \mbox{ for }
n\ge 2.
\end{array}
\right.
\end{equation}
Hence, $(u_n)$ is a sequence of positive numbers and satisfies the
equation
\begin{equation*}
u_{n+1}=u_n+\frac{b(1-1/n)}{\prod_{i=1}^{n}a(1-1/i)},
\end{equation*}
hence
\begin{equation*}
u_n=u_1+\sum_{j=1}^{n-1}\frac{b(1-1/j)}{\prod_{i=1}^{j}a(1-1/i)}.
\end{equation*}
This equality may be rewritten
\begin{equation}\label{expresG}1-G\left(1-\frac{1}{n}\right)=\prod_{i=1}^{n-1}a\left(1-\frac{1}{i}\right)
\left(1-G(0)+\sum_{j=1}^{n-1}\frac{b(1-1/j)}{\prod_{i=1}^{j}a(1-1/i)}\right).
\end{equation}
Using Lemma \ref{LemmaFormeAB}, we easily obtain
\begin{equation}\label{equivproda}\prod_{i=1}^{n}a\left(1-\frac{1}{i}\right)=
\frac{\cste{4}}{n^\alpha}\left(1+\mathcal{O}\left(\frac{1}{n}\right)\right),\quad
\mbox{ with } \cste{4}>0.
\end{equation}
Lemma \ref{LemmaFormeAB} also states that, when $b$ is not
identically $0$ then there exists a unique $k\in\{1,2,\ldots\}$ such
that
\begin{equation*}
\label{devb}b(1-x)=D_k x^k + \mathcal{O}(x^{k+1}), \quad \mbox{ with
} D_k\neq 0.
\end{equation*}
If $b$ is identically $0$, we use the convention $k=+\infty$. In
particular, when $k$ is finite, using (\ref{equivproda}) we deduce
that
\begin{equation}\label{equirap}
\frac{b(1-1/n)}{\prod_{i=1}^{n}a(1-1/i)}=D_k
\cste{4}^{-1}n^{\alpha-k}+\mathcal{O}(n^{\alpha-k-1}).
\end{equation}
Let us now suppose that $k=1$. Combining (\ref{expresG}),
(\ref{equivproda}) and (\ref{equirap}) we find that
$1-G(1-\frac{1}{n})$ converges towards $\frac{D_1}{\alpha}\neq 0$ as
$n$ tends to infinity but this cannot happen because $G$ is
continuous at $1^{-}$ with $G(1)=1$. Thus, we have shown that in
fact
\begin{equation*}
k\geq 2.
\end{equation*}
 We now consider the three cases $\alpha>1$, $\alpha=1$,
$\alpha<1$ separately.
 \vspace*{0.5cm}

\noindent$\boxed{\alpha>1}$

 We have three sub-cases: either
$\alpha>k-1$, or $\alpha<k-1$, or $\alpha = k-1$ with $k\geq 3$.
\begin{itemize}
\item \underline{$\alpha>k-1$}: Making use of (\ref{equirap}), we
have
$$\sum_{j=1}^{n-1}\frac{b(1-1/j)}{\prod_{i=1}^{j}a(1-1/i)}=
\frac{D_k\cste{4}^{-1}}{\alpha-k+1}n^{\alpha-k+1}+\mathcal{O}(1\vee
n^{\alpha-k}).$$ By (\ref{expresG}) and (\ref{equivproda}), we
deduce that
$$1-G\left(1-\frac{1}{n}\right)=
\frac{D_k}{(\alpha-k+1)n^{k-1}}+\mathcal{O}\left(\frac{1}{n^{k\wedge
\alpha}}\right).$$ If $k$ was strictly larger that $2$, we would
have
$$\lim_{n\to\infty} n(1-G(1-1/n)) = 0$$ and therefore
$G'(1)=\mathbb{E}[Z_\infty]=0$ which cannot be true because $Z$ is a
positive random variable which is not equal to zero almost surely.
Thus $k$ must be equal to $2$ and
\begin{equation}\label{espzinfini}1-G\left(1-\frac{1}{n}\right)=
\frac{D_2}{(\alpha-1)n}+\mathcal{O}\left(\frac{1}{n^{2\wedge
\alpha}}\right).\end{equation}

\item \underline{$\alpha<k-1$}:  We prove that this case never happens.
Indeed, in view of ($\ref{equirap}$) we find that, for any
$\varepsilon\in (0,k-1-\alpha)$
\begin{equation}
\frac{b(1-1/n)}{\prod_{i=1}^{n}a(1-1/i)}=\mathcal{O}\left(\frac{1}{n^{1+\varepsilon}}\right)
\end{equation}
(this result also trivially holds when $k=\infty$), thus
\begin{equation*}
\sum_{j=1}^{\infty}\frac{b(1-1/j)}{\prod_{i=1}^{j}a(1-1/i)} <
\infty.
\end{equation*}

Combining this with (\ref{expresG}) and (\ref{equivproda}) we see
that
\begin{equation*}
1-G\left(1-\frac{1}{n}\right)=
\mathcal{O}\left(\frac{1}{n^\alpha}\right).\end{equation*} Since
$\alpha>1$, just as in the previous case, this implies that
$\mathbb{E}[Z_\infty]=0$ which is absurd.

\item \underline{$\alpha=k-1$ and $k\geq 3$}: Again, we prove that
this case is empty. Using (\ref{equirap}), we now get
$$\frac{b(1-1/n)}{\prod_{i=1}^{n}a(1-1/i)}\sim\frac{D_k
\cste{4}^{-1}}{n}.$$
 And, by (\ref{expresG}) and (\ref{equivproda}), we conclude that
$$1-G\left(1-\frac{1}{n}\right)\sim
D_k\frac{ \ln n}{n^{k-1}}.$$ Since $k\geq 3$, we obtain
$\mathbb{E}[Z_\infty]=0$ which is unacceptable.
\end{itemize}
Thus, we have completed the proof of the proposition when $\alpha>1$
and we proved by the way that $k$ must be equal to $2$.

\vspace*{0.5cm}

\noindent$\boxed{\alpha=1}$

We first prove, just as in the previous cases, that $k=2$. Let us
suppose that $k\geq 3$. In view of Lemma \ref{LemmaFormeAB}, for any
$l\geq 3$, we can write the Taylor expansion of $b$ of order $l$
near $1$ in the form
\begin{equation}\label{refpart1}
b(1-x)=D_3x^3+\ldots+D_lx^l+\mathcal{O}(x^{l+1})
\end{equation}
where $D_i\in\R$ for $i\in\{3,4,\ldots,l\}$. Similarly,
\begin{equation*}
a(1-x)=1-x+a_2x^2+\ldots+a_{l}x^{l}+\mathcal{O}(x^{l+1}),
\end{equation*}
from which we deduce that, as $n$ goes to infinity
\begin{equation}\label{refpart2}
\prod_{i=1}^{n}a\left(1-\frac{1}{i}\right)=\frac{a_1'}{n}+
\frac{a_2'}{n^2}+\ldots+\frac{a_{l}'}{n^{l}}+\mathcal{O}
\left(\frac{1}{n^{l+1}}\right)\quad\hbox{with
$a_1'>0$.}
\end{equation}
From (\ref{refpart1}) and (\ref{refpart2}) we also deduce that
\begin{equation*}
\frac{b(1-1/n)}{\prod_{i=1}^{n}a(1-1/i)} =
\frac{d'_2}{n^2}+\ldots+\frac{d'_{l-1}}{n^{l-1}}+\mathcal{O}\left(\frac{1}{n^{l}}\right).
\end{equation*}
Thus,
\begin{equation}\label{refpart3}
\sum_{j=1}^{n-1}\frac{b(1-1/j)}{\prod_{i=1}^{j}a(1-1/i)} = g_0
+\frac{g_1}{n} +
\frac{g_2}{n^2}+\ldots+\frac{g_{l-2}}{n^{l-2}}+\mathcal{O}\left(\frac{1}{n^{l-1}}\right).
\end{equation}
Therefore, in view of (\ref{expresG}), (\ref{refpart2}) and
(\ref{refpart3}), we get
\begin{eqnarray*}
1-G\left(1-\frac{1}{n}\right)&=&\prod_{i=1}^{n-1}a\left(1-\frac{1}{i}\right)
\left(1-G(0)+\sum_{j=1}^{n-1}\frac{b(1-1/j)}{\prod_{i=1}^{j}a(1-1/i)}\right)\\
&=&\frac{\lambda_1}{n}+\frac{\lambda_2}{n^2}+\ldots+
\frac{\lambda_{l-1}}{n^{l-1}}+\mathcal{O}\left(\frac{1}{n^{l}}\right).
\end{eqnarray*}
Comparing with the Taylor expansion of the p.g.f. $G$, we conclude
that $\E(Z_\infty^{l-1}) <\infty$ for all $l$ which contradicts
Proposition \ref{espZM}. Thus, $k=2$ and (\ref{equirap}) yields
$$\frac{b(1-1/n)}{\prod_{i=1}^{n}a(1-1/i)}\sim\frac{D_2
\cste{4}^{-1}}{n} \quad \mbox{ with $D_2\neq 0$.}$$
 And, by (\ref{expresG}) and (\ref{equivproda}), we conclude that
$$1-G\left(1-\frac{1}{n}\right)\sim
D_2\frac{ \ln n}{n},$$ and therefore
\begin{equation}\label{lasted}
\mathbb{E}[Z_\infty]=+\infty.
\end{equation}

\vspace*{0.5cm}

\noindent$\boxed{\alpha<1}$

Since $k\geq 2$, the relation (\ref{equirap}) yields
\begin{equation*}
\sum_{j=1}^{\infty}\frac{b(1-1/j)}{\prod_{i=1}^{j}a(1-1/i)} < \infty
\end{equation*}
(of course, this is trivially true when $k=\infty$). Thus, the
sequence $(u_n)$ defined by (\ref{defun}) converges to a constant
$\cste{5}\geq 0$. Suppose first that $\cste{5}=0$. In this case, $k$
cannot be infinite (because when $k=\infty$, the sequence $(u_n)$ is
constant and then $c_5 = u_1
> 0$). From (\ref{equirap}) we deduce that
$$u_n=-\sum_{j=n}^{\infty}\frac{b(1-1/j)}{\prod_{i=1}^{j}a(1-1/i)}\sim\frac{D_k}{(k-\alpha-1)c_4 n^{k-\alpha-1}},$$
therefore, with the help of (\ref{equivproda}) we get that
$$1-G\left(1-\frac{1}{n}\right)= u_n\prod_{i=1}^{n-1}a\left(1-\frac{1}{i}\right)\sim \frac{D_k}{(k-\alpha-1)n^{k-1}}.$$
Since $k\geq 2$, this implies that $n(1-G(1-1/n))$ converges to a
finite constant and so $\E[Z_\infty]<\infty$. We already notice that
this implies a strict positive speed for the cookie random walk in
the associated cookie environment $(M,\p)$. But (by possibly
extending the value of $M$) we can always construct a cookie
environment $(M,\q)$ such that $\p\leq \q$ and $\alpha(\q)=1$. In
view of (\ref{lasted}), the associated cookie random walk has a zero
speed and this contradicts a monotonicity result of Zerner
(\emph{c.f.} Theorem $17$ of \cite{Zerner05}). Therefore $c_5$
cannot be $0$ and  by (\ref{defun}) and (\ref{equivproda}), we get
that
$$1-G\left(1-\frac{1}{n}\right) = u_n\prod_{i=1}^{n-1}a\left(1-\frac{1}{i}\right) \sim
\frac{\cste{5}\cste{4}}{n^{\alpha}}.$$
\end{proof}

As we already noticed, Theorem \ref{MainTheo1} is now a direct
consequence of the last proposition and Corollary \ref{CorVitesse}.
We also proved that, when the speed is strictly positive, its value
is given by the formula
\begin{equation*}
v=\frac{\alpha-1}{\alpha-1+2D_2}\quad\hbox{where $2D_2=b''(1)>0$.}
\end{equation*}

\begin{rem} In the transient case and when the limiting speed is zero,
Proposition \ref{propG} gives with the help of a classical
Abelian/Tauberian Theorem the asymptotic of the distribution tail of
$Z_\infty$ \emph{i.e.} the distribution tail of the total number of
jumps from $0 $ to $-1$:
\begin{equation}\label{taildist}
\mathbb{P}\left\{Z_\infty>n\right\} \underset{n\to\infty}{\sim} \left\{
\begin{array}{ll}
\frac{c_6}{n^\alpha}&\hbox{if $0<\alpha<1$,}\\
\frac{c_7\ln n }{n}&\hbox{if $\alpha=1$.}
\end{array}\right.
\end{equation}
The functional equation given in Lemma \ref{relationG} for the
p.g.f. of $Z_\infty$ also gives a similar equation for the total
number of returns $R$ to the origin for the cookie random walk.
Indeed, recall that $U^n_0$ (resp. $U^n_1$) stands for the
respective total number of jumps from $0$ to $-1$ (resp. from $1$ to
$0$) before reaching level $n$. Thus, the total number of returns to
the origin before reaching level $n$ is $U^n_0 + U^n_1$ which, under
$\P$ has the same distribution as $Z_n+Z_{n-1}$ under $\mathbb{P}$.
Therefore, we can express the p.g.f. $H$ of the random variable $R$
in term of $G$:
\begin{eqnarray*}
H(s) &=& \mathbb{E}\left[s^{Z_\infty}\mathbb{E}_{Z_\infty}\left[s^{Z_\infty}\right]\right]\\
&=& \frac{1}{a(s)}G\left(\frac{s}{2-s}\right) +
\sum_{k=0}^{M-2}G^{(k)}(0)s^k\left(\mathbb{E}\left[s^{A_k}\right] - \frac{1}{a(s)(2-s)^k}\right).
\end{eqnarray*}
In particular, Proposition \ref{propG} also holds for $H$ and the
tail distribution of the total number of returns to the origin when
$\alpha\leq 1$ has the same form as in (\ref{taildist}).
\end{rem}

\begin{rem}In the particular case $M=2$ (there are at most $2$ cookies per
site), the only unknown in the definition of the function $b$ is
$G(0)$. Since we know that $b'(1)=0$ (\emph{c.f.} the beginning of
the proof of Proposition \ref{propG}) we can therefore explicitly
calculate $G(0)$, that is the probability that the cookie random
walk never jumps from $0$ to $1$ which is also the probability that
the cookie random walk never hits $-1$. According to the previous
remark, we can also calculate the probability that the cookie random
walk never returns to $0$. Hence, we recover Theorem $18$ of
\cite{Zerner05} in the case of a deterministic cookie environment.
\end{rem}

\section{Continuity of the speed and differentiability at the critical point}
The aim of this section is to prove Theorem \ref{MainTheo2}. Recall that
\begin{equation*}
v(M,\p)=\left\{
\begin{array}{ll}
0&\hbox{ if $\alpha(M,p)\leq 1$,}\\
\frac{\alpha-1}{\alpha-1+b''(1)}&\hbox{ if $\alpha(M,p) > 1$,}
\end{array}\right.
\end{equation*}
where $b''(1)$ stands for the second derivative at point $1$ of the
function $b$ defined in Lemma \ref{relationG}:
\begin{equation*}
b(s)=1-\frac{1}{(2-s)^{M-1}\E\left[s^{A_{M-1}}\right]}+
 \sum_{k=0}^{M-2}\mathbb{P}\{Z_\infty=k\}
 \left(\frac{\E\left[s^{A_{k}}\right]}{(2-s)^{M-1}\E\left[s^{A_{M-1}}\right]}-\frac{1}{(2-s)^k}
 \right).
\end{equation*}
Furthermore, we also proved in Proposition \ref{propG} that, when
$\alpha(M,\p)=1$, then $b''(1)$ is strictly positive. Hence, in
order to prove Theorem \ref{MainTheo2}, we just need to show that
$b''(1) = b''_{(M,\p)}(1)$ is a continuous function of $\p$ in
$\Omega^u_{M}$. It is also clear from the definition of the random
variables $A_k$ that the functions
\begin{equation*}
\p \to
\left(\E_{(M,\p)}\left[s^{A_k}\right]\right)^{(i)}(1)\quad\hbox{(\emph{i.e.}
the $i^{\hbox{\tiny{th}}}$ derivative at point $1$)}
\end{equation*}
are continuous in $\p$ in $\Omega^u_{M}$ for all $k\geq 0$ and all
$i\geq 0$ (it is a rational function in $p_1,\ldots,p_M$).
Therefore, it simply remains to prove that, for any $k\geq 0$, the
function
\begin{equation*}
\p\to\mathbb{P}_{(M,\p)}\left\{Z_\infty = k\right\}
\end{equation*}
is continuous in $\Omega^u_M$. The following lemma is based on the
monotonicity of the hitting times of a cookie random walk with
respect to the environment.

\begin{duge}
Let $(M,\p)$ be a cookie environment such that $\alpha(M,\p)>0$.
Then there exist $\varepsilon>0$ and $f:\N\mapsto \R_+$ with
$\lim_{n\rightarrow +\infty}f(n)=0$ such that
$$\forall \q\in B(\p,\varepsilon),\forall j\in \N,\forall n\in\N,
\;|\mathbb{P}_{(M,\q)}\left\{Z_\infty=j\right\}-\mathbb{P}_{(M,\q)}\left\{Z_n=j\right\}|\le
f(n),$$ where
$$
B(\p,\varepsilon)=\Big\{\q=(q_1,\ldots,q_M),\; \frac{1}{2}\le q_i
<1,\; \alpha(M,\q)>0 \mbox{ and } \sum_{i=1}^\infty|p_i- q_i|\le
\varepsilon\Big\}.
$$
\end{duge}

\begin{proof}
Let us fix $(M,\p)$ with $\alpha(M,\p)>0$. For $\varepsilon>0$,
define the vector
$\p^{\varepsilon}=(p_1^{\varepsilon},\ldots,p_M^{\varepsilon})$ by
$p_i^{\varepsilon}=\max(\frac{1}{2},p_i-\varepsilon)$. We can choose
$\varepsilon >0$ such that $\alpha(M,\p^{\varepsilon})>0$. Then, for
all $\q \in B(\p,\varepsilon)$, we have
\begin{equation}\label{onemort}
\p^{\varepsilon}\leq \q
\end{equation}
(where $\leq$ denotes the canonical partial order on $\R^M$). Let
now pick $\q\in B(\p,\varepsilon)$, $j\in \N$ and $n\in \N$. Recall
that $U^\infty_0$ denotes the total number of jump of the cookie
random walk from $0$ to $-1$ and
$$\mathbb{P}_{(M,\q)}\{Z_\infty=j\}=\P_{(M,\q)}\{U_0^{\infty}=j\}=
\P_{(M,\q)}\{X \mbox{ jumps }j \mbox{ times from 0 to -1}\},$$
and
\begin{multline*}
\mathbb{P}_{(M,\q)}\{Z_n=j\}=\P_{(M,\q)}\{U_0^{n}=j\}\\=\P_{(M,\q)}\{X
\mbox{ jumps }j \mbox{ times from 0 to -1 before reaching } n\}.
\end{multline*}
Hence
\begin{eqnarray}
\nonumber|\mathbb{P}_{(M,\q)}\{Z_\infty=j\}-\mathbb{P}_{(M,\q)}\{Z_n=j\}|
&=&|\P_{(M,\q)}\{U_0^{\infty}=j\}-\P_{(M,\q)}\{U_0^{n}=j\}|\\
\nonumber&\le& \P_{(M,\q)}\{U_0^{n}\neq U_0^\infty\}\\
&= & \P_{(M,\q)}\{A\},\label{oncerr}
\end{eqnarray}
where $A$ is the event "$X \mbox{ visits } -1 \mbox{ at least once
after reaching level } n$". Recall the notation $\omega =
\omega(i,x)_{i\geq 1,x\in\Z}$ for a general cookie environment given
in the introduction. Let now $\omega_{X,n}$ denote the (random)
cookie-environment obtained when the cookie random walk $X$ hits
level $n$ for the first time and shifted by $n$, i.e. for all $x\in
\Z$ and $i\ge 1$, if the initial cookie environment is $\omega$,
then
\begin{equation*}
\omega_{X,n}(i,x)=\omega(j,x+n) \quad \mbox{ where }
j=i+\sharp\{0\le k< T_n, X_k=x+n\}.
\end{equation*}
With this notation we have
$$\P_{(M,\q)}\left\{A\right\}=\E_{(M,\q)}
\left[\P_{\omega_{X,n}}\{X \mbox{ visits $-(n+1)$ at least
once}\}\right].$$ Besides, $X$ has not eaten any cookie at the sites
$x\geq n$ before time $T_n$. Thus, the environment $\omega_{X,n}$
satisfies $\P_{(M,\q)}$-almost surely
$$\omega_{X,n}(i,x)=q_i, \quad \mbox{ for all } x\ge 0 \mbox{ and } i\geq
1 \hbox{ (with the convention $q_i=\frac{1}{2}$ for $i>M$).}$$

Hence, in view of (\ref{onemort}), the random cookie environment
$\omega_{X,n}$ is $\P_{(M,\q)}$-almost surely larger (for the
canonical partial order) than the deterministic environment
$\omega_{\p^{\varepsilon}}$ defined by
\begin{equation*} \left\{
\begin{array}{l}
\omega_{\p^\varepsilon}(i,x)=\frac{1}{2}, \quad \mbox{ for all } x<
0 \mbox{ and }
i\ge 1, \\
\omega_{\p^\varepsilon}(i,x)=p^{\varepsilon}_i, \quad \mbox{ for all
} x\ge 0 \mbox{ and } i\ge 1 \hbox{ (with the convention
$p^{\varepsilon}_i=\frac{1}{2}$ for $i\geq M$).}
\end{array}\right.
\end{equation*}
Thus, Lemma $15$ of \cite{Zerner05} yields
\begin{multline*}
\P_{\omega_{X,n}}\{X \mbox{ visits } -(n+1)\mbox{ at least once}\}
\\\leq \P_{\omega_{\p^{\varepsilon}}}\{X \mbox{ visits } -(n+1)\mbox{
at least once}\}\quad \P_{(M,\q)}-a.s.
\end{multline*}
 In view of (\ref{oncerr}) we
deduce that
$$|\mathbb{P}_{(M,\q)}\{Z_\infty=j\}-\mathbb{P}_{(M,\q)}\{Z_n=j\}|\le f(n),$$
where $f(n)=\P_{\omega_{\p^{\varepsilon}}}\{X \mbox{ visits $-(n+1)$
at least once}\}$ does not depend of $\q$. It remains to prove that
$f(n)$ tends to $0$ as $n$ goes to infinity. Let us first notice
that
$$\P_{\omega_{\p^{\varepsilon}}}\{\forall n\ge 0,\; X_n\ge
0\}=\P_{(M,\p^{\varepsilon})}\{\forall n\ge 0,\; X_n\ge 0\},$$ since
these probabilities depend only on the environments on the half line
$[0,+\infty)$. Recall also that the cookie random walk in the
environment $(M,\p^{\varepsilon})$ is transient (we have chosen
$\varepsilon$ such that $\alpha(M,\p^{\varepsilon})>0$), thus
$$\P_{(M,\p^{\varepsilon})}\{\forall n\ge 0,\; X_n\ge
0\}=\P_{(M,\p^{\varepsilon})}\{U_0^\infty=0\}=\mathbb{P}_{(M,\p^{\varepsilon})}\{Z_\infty=0\}>0.$$
Hence
$$\P_{\omega_{\p^{\varepsilon}}}\{\forall n\ge 0,\; X_n\ge
0\}>0,$$
which implies
$$
\P_{\omega_{\p^{\varepsilon}}}\{\hbox{$X_n=0$ infinitely often}\}<1,
$$
and a $0-1$ law (\emph{c.f.} Proposition $5$ of \cite{Zerner05}) yields
$$
\P_{\omega_{\p^{\varepsilon}}}\{\hbox{$X_n=0$ infinitely
often}\}=\P_{\omega_{\p^{\varepsilon}}}\{\hbox{$X_n\leq 0$
infinitely often}\}=0.
$$
Therefore, $\lim_{n\to\infty}f(n)=0$.
\end{proof}

Recall that the transition probabilities of the Markov chain $Z$ are
given by the law of the random variables $A_k$:
\begin{equation*}
\mathbb{P}_{(M,\p)}\left\{Z_{n+1} = j\;|\;Z_n=i\right\} =
\P_{(M,\p)}\left\{A_i=j\right\}.
\end{equation*}
It is therefore clear that for each fixed $n$ and each $k$, the
function $\p\to\mathbb{P}_{(M,\p)}\left\{Z_n = k\right\}$ is
continuous in $\p$ in $\Omega^u_{M}$. In view of the previous lemma,
we conclude that for each $k$ the function
$\p\to\mathbb{P}_{(M,\p)}\left\{Z_\infty = k\right\}$ is also
continuous in $\p$ in $\Omega^u_{M}$ and this completes the proof of
Theorem \ref{MainTheo2}.

\begin{ack}The authors would like to thank Yueyun Hu for all his precious advices.
\end{ack}

\nocite{Antal05}
\nocite{Zerner06}
\bibliographystyle{plain}
\bibliography{bibliocookie}

\end{document}